\documentclass[12pt]{amsart}
\usepackage[all]{xy}

\usepackage{amsmath, amsthm, amssymb,latexsym}

\newcommand{\N}{\mathbb{N}}

\newcommand{\R}{\mathbb{R}}
\newcommand{\C}{\mathbb{C}}

\newcommand{\mR}{\mathcal{R}}

\newcommand{\norm}[1]{\Vert#1\Vert}
\newcommand{\bignorm}[1]{\bigl\Vert#1\bigr\Vert}
\newcommand{\Bignorm}[1]{\Bigl\Vert#1\Bigr\Vert}
\newcommand{\biggnorm}[1]{\biggl\Vert#1\biggl\Vert}

\newcommand{\HI}{H^\infty}

\DeclareMathOperator{\bi}{\!^{**}}

\DeclareMathOperator{\BI}{\mathcal{B}^\infty}

\DeclareMathOperator{\Rad}{Rad}

\newtheorem{thm}{Theorem}[section]
\newtheorem{lem}[thm]{Lemma}
\newtheorem{cor}[thm]{Corollary}
\newtheorem{prop}[thm]{Proposition}
\newtheorem{defi}[thm]{Definition}
\theoremstyle{remark}
\newtheorem{rem}[thm]{\bf Remark}
\numberwithin{equation}{section}

\begin{document}

\title[]{Tensor extension properties of $C(K)$-representations
and applications to unconditionality}

\author{Christoph Kriegler, Christian Le Merdy}
\address{Institut f\"ur Analysis\\ Kaiserstrasse 89 \\ 76133 Karlsruhe\\
Germany}

\address{Laboratoire de Math\'ematiques\\ Universit\'e de
Franche-Comt\'e
\\ 25030 Besan\c con Cedex\\ France}
\email{christoph.kriegler@univ-fcomte.fr}

\address{Laboratoire de Math\'ematiques\\ Universit\'e de  Franche-Comt\'e
\\ 25030 Besan\c con Cedex\\ France}
\email{clemerdy@univ-fcomte.fr}

\date{\today}

\thanks{The first author is supported
by the Karlsruhe House of Young Scientists and the Franco-German
University DFH-UFA, the second author is supported by the research
program ANR-06-BLAN-0015.}

\begin{abstract} Let $K$ be any compact set. The $C^*$-algebra $C(K)$ is nuclear and
any bounded homomorphism from $C(K)$ into $B(H)$, the algebra of
all bounded operators on some Hilbert space $H$, is automatically
completely bounded. We prove extensions of these results to the
Banach space setting, using the key concept of $R$-boundedness.
Then we apply these results to operators with a uniformly bounded
$\HI$-calculus, as well as to unconditionality on $L^p$. We show
that any unconditional basis on $L^p$ `is'  an unconditional basis
on $L^2$ after an appropriate change of density.
\end{abstract}

\maketitle


\bigskip

\section{Introduction}\label{Sec Intro}
Throughout the paper, we let $K$ be a compact set and we let
$C(K)$ be the algebra of all continuous functions $f\colon
K\to\C$, equipped with the supremum norm. A representation of
$C(K)$ on some Banach space $X$ is a bounded unital homomorphism
$u\colon C(K)\to B(X)$ into the algebra $B(X)$ of all bounded
operators on $X$. Such representations appear naturally and play a
major role in several fields of Operator Theory, including
functional calculi, spectral theory and spectral measures, or
classification of $C^*$-algebras. Several recent papers, in
particular \cite{KaW, KW, DG, dPR}, have emphasized the rich and
fruitful interplays between the notion of $R$-boundedness,
uncondionality and various functional calculi. The aim of this
paper is to establish new properties of $C(K)$-representations
involving $R$-boundedness, and to give applications to
$\HI$-calculus (in the sense of \cite{CDMY,KaW}) and to
unconditionality in $L^p$-spaces.

We recall the definition of $R$-boundedness (see \cite{BG,CPSW}).
Let $(\epsilon_k)_{k\geq 1}$ be a sequence of independent
Rademacher variables on some probability space $\Omega_0$. That
is, the $\epsilon_k$'s take values in $\{-1,1\}$ and ${\rm
Prob}(\{\epsilon_k =1\}) = {\rm Prob}(\{\epsilon_k
=-1\})=\frac{1}{2}$. For any Banach space $X$, we let
$\Rad(X)\subset L^2(\Omega_0;X)$  be the closure of ${\rm
Span}\{\epsilon_k\otimes x\, :\, k\geq 1,\ x\in X\}$ in
$L^2(\Omega_{0};X)$. Thus for any finite family $x_1,\ldots,x_n$
in $X$, we have
$$
\Bignorm{\sum_k \epsilon_k\otimes x_k}_{\Rad(X)} \,=\,
\Bigr(\int_{\Omega_0}\Bignorm{\sum_k \epsilon_k(\lambda)\,
x_k}_{X}^{2}\,d\lambda\,\Bigr)^{\frac{1}{2}}.
$$
By definition, a set $\tau\subset B(X)$ is $R$-bounded if there is
a constant $C\geq 0$ such that for any finite families
$T_1,\ldots, T_n$ in $\tau$, and $x_1,\ldots,x_n$ in $X$, we have
$$
\Bignorm{\sum_k \epsilon_k\otimes T_k x_k}_{\Rad(X)}\,\leq\, C\,
\Bignorm{\sum_k \epsilon_k\otimes x_k}_{\Rad(X)}.
$$
In this case, we let $R(\tau)$ denote the smallest possible $C$.
It is called the $R$-bound of $\tau$. By convention, we write
$R(\tau)=\infty$ if $\tau$ is not $R$-bounded.

It will be convenient to let $\Rad_n(X)$ denote the subspace of
$\Rad(X)$ of all finite sums $\sum_{k=1}^{n}\epsilon_k\otimes
x_k$. If $X=H$ is a Hilbert space, then  $\Rad_n(H)=\ell^2_n(H)$
isometrically and all bounded subsets of $B(H)$ are automatically
$R$-bounded. On the opposite, if $X$ is not isomorphic to a
Hilbert space, then $B(X)$ contains bounded subsets which are not
$R$-bounded \cite[Prop. 1.13]{AB}.

\bigskip
As a guise of motivation for this paper, we recall two well-known
properties of $C(K)$-representations on Hilbert space $H$. First,
any bounded homomorphism $u\colon C(K)\to B(H)$ is completely
bounded, with $\norm{u}_{cb}\leq\norm{u}^2$. This means that for
any integer $n\geq 1$, the tensor extension $I_{M_n}\otimes
u\colon M_n(C(K))\to M_n(B(H))$ satisfies $\norm{I_{M_n}\otimes
u}\leq\norm{u}^2$, if $M_n(C(K))$ and $M_n(B(H))$ are equipped
with their natural $C^*$-algebra norms. This implies that any
bounded homomorphism $u\colon C(K)\to B(H)$ is similar to a
$*$-representation, a result going back at least to \cite{Bu}. We
refer to \cite{Pa,Pis2} and the references therein for some
information on completely bounded maps and similarity properties.

Second, let $u\colon C(K)\to B(H)$ be a bounded homomorphism. Then
for any $b_1,\cdots,b_n$ lying in the commutant of the range of
$u$ and for any $f_1,\ldots, f_n$ in $C(K)$, we have
\begin{equation}\label{Nuclear}
\Bignorm{\sum_k u(f_k)b_k}\,\leq\,\norm{u}^{2}\sup_{t\in K}
\Bignorm{\sum_k f_k(t)b_k}.
\end{equation}
This property is essentially a rephrasing of the fact that $C(K)$
is a nuclear $C^*$-algebra. More precisely, nuclearity means that
the above property holds true for $*$-representations (see e.g.
\cite[Chap. 11]{KR} or \cite[Chap. 12]{Pa}), and its extension to
arbitrary bounded homomorphisms easily follows from the similarity
property mentioned above (see \cite{LM2} for more explanations and
developments).

\smallskip
Now let $X$ be a Banach space and let $u\colon C(K)\to B(X)$ be a
bounded homomorphism. In Section 2, we will show the following
analog of (\ref{Nuclear}):
\begin{equation}\label{Nuclear2}
\Bignorm{\sum_k u(f_k)b_k}\,\leq\,\norm{u}^{2}
R\Bigl(\Bigl\{\sum_k f_k(t) b_k: t\in K\Bigr\}\Bigr),
\end{equation}
provided that the $b_k$'s commute with the range of
$u$.

\smallskip
Section 3 is devoted to sectorial operators $A$ which
have a uniformly bounded $\HI$-calculus, in the sense that they
satisfy an estimate
\begin{equation}\label{Estimate}
\norm{f(A)}\leq C\sup_{t>0}\vert f(t)\vert
\end{equation}
for any bounded analytic function on a sector $\Sigma_\theta$
surrounding $(0,\infty)$. Such operators turn out to have a
natural $C(K)$-functional calculus. Applying (\ref{Nuclear2}) to
the resulting representation $u\colon C(K)\to B(X)$, we show that
(\ref{Estimate}) can be automatically extended to operator valued
analytic functions $f$ taking their values in the commutant of
$A$. This is an analog of a remarkable result of Kalton-Weis
\cite[Thm. 4.4]{KW} saying that if an operator $A$ has a bounded
$\HI$-calculus and $f$ is an operator valued analytic function
taking its values in an $R$-bounded subset of the commutant of
$A$, then the operator $f(A)$ arising from `generalized
$\HI$-calculus', is bounded.

\smallskip
In Section 4, we introduce matricially $R$-bounded maps $C(K)\to
B(X)$, a natural analog of completely bounded maps in the Banach
space setting. We show that if $X$ has property $(\alpha)$, then
any bounded homomorphism $C(K)\to B(X)$ is automatically
matricially $R$-bounded. This extends both the Hilbert space
result mentioned above, and a result of De Pagter-Ricker
\cite[Cor. 2.19]{dPR} saying that any bounded homomorphism
$C(K)\to B(X)$ maps the unit ball of $C(K)$ into an $R$-bounded
set, provided that $X$ has property $(\alpha)$.

\smallskip
In Section 5, we give an application of matricial $R$-boundedness
to the case when $X=L^p$. A classical result of Johnson-Jones
\cite{JJ} asserts that any bounded operator $T\colon L^p\to L^p$
acts, after an appropriate change of density, as a  bounded
operator on $L^2$. We show versions of this theorem for bases
(more generally, for FDD's). Indeed we show that any unconditional
basis (resp. any $R$-basis) on $L^p$ becomes an unconditional
basis (resp. a Schauder basis) on $L^2$ after an appropriate
change of density. These results rely on Simard's extensions of
the Johnson-Jones Theorem established in \cite{Sim}.

\bigskip
We end this introduction with a few preliminaries and notation.
For any Banach space $Z$, we let $C(K;Z)$ denote the space of all
continuous functions $f\colon K\to Z$, equipped with the supremum
norm
$$
\norm{f}_{\infty}=\sup\{\norm{f(t)}_Z\, :\, t\in K\}.
$$
We may regard $C(K)\otimes Z$ as a subspace of $C(K;Z)$, by
identifying $\sum_k f_k\otimes z_k$ with the function
$t\mapsto\sum_k f_k(t)z_k$, for any finite families $(f_k)_k$ in
$C(K)$ and $(z_k)_k$ in $Z$. Moreover, $C(K)\otimes Z$ is dense in
$C(K;Z)$. Note that for any integer $n\geq 1$, $C(K;M_n)$
coincides with the $C^*$-algebra $M_n(C(K))$ mentioned above.

We will need the so-called `contraction principle', which says
that for any $x_1,\ldots, x_n$ in a Banach space $X$ and any
$\alpha_1,\ldots,\alpha_n$ in $\C$, we have
\begin{equation}\label{CP}
\Bignorm{\sum_k \epsilon_k\otimes \alpha_k x_k}_{\Rad(X)}\leq
2\,\sup_k\vert \alpha_k\vert\, \Bignorm{\sum_k \epsilon_k\otimes
x_k}_{\Rad(X)}.
\end{equation}

We also recall that any commutative $C^*$-algebra is a
$C(K)$-space (see e.g. \cite[Chap. 4]{KR}). Thus our results
concerning $C(K)$-representations apply as well to all these
algebras. For example we will apply them to $\ell^\infty$ in
Section 5.

We let $I_X$ denote the identity mapping on a Banach space $X$,
and we let $\chi_B$ denote the indicator function of a set $B$. If
$X$ is a dual Banach space, we let $w^*B(X)\subset B(X)$ be the
subspace of all $w^*$-continuous operators on $X$.

\medskip

\section{The extension theorem}\label{Sec Main}

Let $X$ be an arbitrary Banach space. For any compact set $K$ and
any bounded homomorphism $u\colon C(K)\to B(X)$, we let
$$
E_u\,=\,\bigl\{ b\in B(X)\, :\, bu(f)=u(f)b\,\ \hbox{for any}\
f\in C(K)\bigr\}
$$
denote the commutant of the range of $u$.

Our main purpose in this section is to prove (\ref{Nuclear}). We
start with the case when $C(K)$ is finite dimensional.

\begin{prop}\label{Prop Main thm for l infty n->B(X)}
Let $N\geq 1$ and let $u\colon \ell^\infty_N\to B(X)$ be a bounded
homomorphism. Let $(e_1,\ldots,e_N)$ be the canonical basis of
$\ell^\infty_N$ and set $p_i = u(e_i)$, $\,i=1,\ldots,N.$ Then for
any $b_1,\ldots,b_N\in E_u$, we have
$$
\Bignorm{\sum_{i=1}^N p_i b_i} \leq \|u\|^2 R(\{b_1,\ldots,b_N\}).
$$
\end{prop}

\begin{proof}
Since $u$ is multiplicative, each $p_i$ is a projection and $p_i
p_j = 0$ when $i\neq j.$ Hence for any choice of signs
$(\alpha_1,\ldots,\alpha_N)\in \{-1,1\}^N,$ we have
$$
\sum_{i=1}^N p_i b_i \, = \, \sum_{i,j=1}^N \alpha_i \alpha_j p_i
p_j b_j.
$$
Furthermore,
$$
\Bignorm{\sum_i \alpha_i p_i} =
\|u(\alpha_1,\ldots,\alpha_N)\|\leq \|u\|\,
\|(\alpha_1,\ldots,\alpha_N)\|_{\ell^\infty_N} = \|u\|.
$$
Therefore for any $x\in X$, we have the following chain of
inequalities which prove the desired estimate:
\begin{align}
\Bignorm{\sum_i p_i b_i x }^{2} & = \, \int_{\Omega_0}
\Bignorm{\sum_i \epsilon_i(\lambda) p_i \sum_j \epsilon_j(\lambda)
p_j b_j x}^{2}\, d\lambda \nonumber\\ & \leq \, \int_{\Omega_0}
\Bignorm{ \sum_i \epsilon_i(\lambda) p_i}^{2}\, \Bignorm{\sum_j
\epsilon_j(\lambda) p_j b_j x}^{2}\, d\lambda \nonumber\\ & \leq
\|u\|^{2}\,\int_{\Omega_0} \Bignorm{\sum_j \epsilon_j(\lambda) b_j
p_j x}^{2}\, d\lambda \nonumber\\ & \leq \|u\|^{2}\,
R(\{b_1,\ldots,b_N\})^{2} \,\int_{\Omega_0} \Bignorm{\sum_j
\epsilon_j(\lambda) p_j x}^{2}\, d\lambda \nonumber\\ & \leq
\|u\|^{4} R(\{b_1,\ldots,b_N\})^{2} \|x\|^{2}. \nonumber
\end{align}
\end{proof}

The study of infinite dimensional $C(K)$-spaces requires the use
of second duals and $w^*$-topologies. We recall a few well-known
facts that will be used later on in this paper.  According to the
Riesz representation Theorem, the dual space $C(K)^*$ can be
naturally identified with the space $M(K)$ of Radon measures on
$K$. Next, the second dual space $C(K)^{**}$ is a commutative
$C^*$-algebra for the so-called Arens product. This product
extends the product on $C(K)$ and is separately $w^*$-continuous,
which means that for any $\xi\in C(K)^{**}$, the two linear maps
$$
\nu \in C(K)^{**}\longmapsto \nu\xi \in
C(K)^{**}\qquad\hbox{and}\qquad \nu \in C(K)^{**}\longmapsto
\xi\nu \in C(K)^{**}
$$
are $w^*$-continuous.

Let
$$
\BI(K)\, =\, \{f\colon K\to \C\,\vert\, f\text{ bounded, Borel
measurable}\},
$$
equipped with the sup norm. According to the duality pairing
$$
\langle f,\mu\rangle\,=\, \int_K f(t)\, d\mu(t),\qquad \mu\in
M(K),\ f\in \BI(K),
$$
one can regard $\BI(K)$ as a closed subspace of $C(K)^{**}$.
Moreover the restriction of the Arens product to $\BI(K)$
coincides with the pointwise product. Thus we have natural
$C^*$-algebra inclusions
\begin{equation}\label{seconddual}
C(K)\subset \BI(K)\subset C(K)^{**}.
\end{equation}
See e.g. \cite[pp. 366-367]{Da} and \cite[Sec. 9]{C} for further
details.

Let $\widehat{\otimes}$ denote the projective tensor product on
Banach spaces. We recall that for any two Banach spaces $Y_1,Y_2$,
we have a natural identification
$$
(Y_1\widehat{\otimes } Y_2)^* \simeq B(Y_2,Y_1^{*}),
$$
see e.g. \cite[VIII.2]{DU}. This implies that when $X$ is a dual
Banach space, $X=(X_*)^*$ say, then
$B(X)=(X_*\widehat{\otimes}X)^*$ is a dual space. The next two
lemmas are elementary.

\begin{lem}\label{Lem Product on B(X*) separately w*-cont}
Let $X=(X_*)^*$ be a dual space, let $S\in B(X)$, and let
$R_S,L_S\colon B(X)\to B(X)$ be the right and left multiplication
operators defined by $R_S(T)= TS$ and $L_S(T)= ST$. Then $R_S$ is
$w^*$-continuous whereas $L_S$ is $w^*$-continuous if (and only
if) $S$ is $w^*$-continuous.
\end{lem}

\begin{proof}
The tensor product mapping $I_{X_*}\otimes S$  on $X_{*}\otimes X$
uniquely extends to a bounded map $r_S\colon
X_*\widehat{\otimes}X\to X_*\widehat{\otimes}X$, and we have
$R_S=r_S^*$. Thus $R_S$ is $w^*$-continuous. Likewise, if $S$ is
$w^*$-continuous and if we let $S_*\colon X_*\to X_*$ be its
pre-adjoint map, the tensor product mapping $S_*\otimes I_{X}$  on
$X_{*}\otimes X$ extends to a bounded map $l_S\colon
X_*\widehat{\otimes}X\to X_*\widehat{\otimes}X$, and $L_S=l_S^*$.
Thus $L_S$ is $w^*$-continuous. The `only if' part (that we will
not use) is left to the reader.
\end{proof}

\begin{lem}\label{Lem Extension C(K)->B(X) to C(K)**->B(X)}
Let $u\colon C(K)\to B(X)$ be a bounded map. Suppose that $X$ is a
dual space. Then there exists a (necessarily unique)
$w^*$-continuous linear mapping  $\widetilde{u}\colon C(K)^{**}\to
B(X)$ whose restriction to $C(K)$ coincides with $u$. Moreover
$\norm{\widetilde{u}}=\norm{u}$.

If further $u$ is a homomorphism, and $u$ is valued in $w^*B(X)$,
then $\widetilde{u}$ is a homomorphism as well.
\end{lem}

\begin{proof}
Let $j\colon (X_*\widehat{\otimes} X) \hookrightarrow
(X_*\widehat{\otimes} X)\bi$ be the canonical injection and
consider its adjoint $p=j^*\colon B(X)\bi\to B(X)$. Then set
$$
\widetilde{u} = p \circ u\bi \colon C(K)\bi \longrightarrow B(X).
$$
By construction, $\widetilde{u}$ is $w^*$-continuous and extends
$u$. The equality $\norm{\widetilde{u}}=\norm{u}$ is clear.

Assume now that $u$ is a homomorphism and that   $u$ is valued in
$w^*B(X)$. Let $\nu,\xi\in C(K)\bi$ and let $(f_\alpha)_\alpha$
and $(g_\beta)_\beta$ be bounded nets in $C(K)$ $w^*$-converging
to $\nu$ and $\xi$ respectively. By both parts of Lemma \ref{Lem
Product on B(X*) separately w*-cont}, we have the following
equalities, where limits are taken in the $w^*$-topology of either
$C(K)^{**}$ or $B(X)$:
$$
\widetilde{u}(\nu\xi)  = \widetilde{u}(\lim_\alpha \lim_\beta
f_\alpha g_\beta) = \lim_\alpha \lim_\beta u(f_\alpha g_\beta) =
\lim_\alpha \lim_\beta u(f_\alpha) u(g_\beta) = \lim_\alpha
u(f_\alpha) \widetilde{u}(\xi)   =
\widetilde{u}(\nu)\widetilde{u}(\xi).
$$
\end{proof}

We refer e.g. to \cite[Lem. 2.4]{HKK} for the following fact.

\begin{lem}\label{Lem R(tau) = R(tau**)}
Consider $\tau \subset B(X)$ and set $\tau\bi = \{T\bi\, :\, T\in
\tau\}\subset B(X\bi).$ Then $\tau$ is $R$-bounded if and only if
$\tau^{**}$ is $R$-bounded, and in this case,
$$
R(\tau) = R(\tau\bi).
$$
\end{lem}

For any $F\in C(K;B(X))$, we set
$$
R(F)=R\bigl(\{F(t)\, :\, t\in K\}\bigr).
$$
Note that $R(F)$ may be infinite. If $F=\sum_k f_k\otimes b_k$
belongs to the algebraic tensor product $C(K)\otimes B(X)$, we set
$$
\Bignorm{\sum_k f_k\otimes b_k}_R = R(F) = R\Bigl(\Bigl\{\sum_k
f_k(t) b_k\, :\, t\in K\Bigr\}\Bigr).
$$
Note that by (\ref{CP}), we have
\begin{equation}\label{R}
\norm{f\otimes b}_R\leq 2\norm{f}_\infty\norm{b},\qquad f\in
C(K),\ b\in B(X).
\end{equation}
From this it is easy to check that $\norm{\ }_R$ is finite and is
a norm on $C(K)\otimes B(X)$.

Whenever $E\subset B(X)$ is a closed subspace, we let
$$
C(K)\mathop{\otimes}\limits^R E
$$
denote the completion of $C(K)\otimes E$ for the norm $\norm{\
}_R$.

\begin{rem}\label{Embed} Since the $R$-bound of a set is greater than its
uniform bound, we have $\norm{\ }_\infty\leq\norm{\ }_R$ on
$C(K)\otimes B(X)$. Hence the canonical embedding of $C(K)\otimes
B(X)$ into $C(K;B(X))$ uniquely extends to a contraction
$$
J\colon C(K)\mathop{\otimes}\limits^R B(X)\longrightarrow
C(K;B(X)).
$$
Moreover $J$ is 1-1 and for any $\varphi\in
C(K)\mathop{\otimes}\limits^R B(X)$, we have
$R(J(\varphi))=\norm{\varphi}_R$. Indeed, let $(F_n)_{n\geq 1}$ be
a sequence in $C(K)\otimes B(X)$ such that
$\norm{F_n-\varphi}_R\to 0$ and let $F=J(\varphi)$. Then
$\norm{F_n}_R \to\norm{\varphi}_R$ and $\norm{F_n-F}_\infty\to 0$.
According to the definition of the $R$-bound, the latter property
implies that $\norm{F_n}_R\to\norm{F}_R$, which yields the result.
\end{rem}

\begin{thm}\label{Thm Main}
Let $u\colon C(K)\to B(X)$ be a bounded homomorphism.

\begin{itemize}
\item [(1)] For any finite families $(f_k)_k$ in $C(K)$ and
$(b_k)_k$ in $E_u$, we have
$$
\Bignorm{\sum_k u(f_k)b_k}\,\leq\,\norm{u}^{2}\,\Bignorm{\sum_k
f_k\otimes b_k}_R.
$$
\item [(2)] There is a (necessarily unique) bounded linear map
$$
\widehat{u}\colon  C(K)\mathop{\otimes}\limits^R
E_u\longrightarrow B(X)
$$
such that $\widehat{u}(f\otimes b)=u(f)b$ for any $f\in C(K)$ and
any $b\in E_u$. Moreover, $\norm{\widehat{u}}\leq \norm{u}^{2}$.
\end{itemize}
\end{thm}

\begin{proof}
Part (2) clearly follows from part (1). To prove (1) we introduce
$$
w\colon C(K) \longrightarrow B(X\bi),\qquad  w(f) = [u(f)]\bi.
$$
Then $w$ is a bounded  homomorphism and $\|w\| = \|u\|$. We let
$\widetilde{w}\colon C(K)^{**}\to B(X^{**})$ be its
$w^*$-continuous extension given by Lemma \ref{Lem Extension
C(K)->B(X) to C(K)**->B(X)}. Note that $w$ is valued in
$w^*B(X^{**})$, so $\widetilde{w}$ is a homomorphism. We claim
that
$$
\{b\bi\, :\,b\in E_u\} \subset E_{\widetilde{w}}.
$$
Indeed, let $b\in E_u$. Then for all $f\in C(K),$ we have
$$
b\bi w(f) = (bu(f))\bi = (u(f)b)\bi = w(f)b\bi.
$$
Next for any $\nu\in C(K)^{**}$, let $(f_\alpha)_\alpha$ be a
bounded net in $C(K)$ which converges to $\nu$ in the
$w^*$-topology. Then by Lemma \ref{Lem Product on B(X*) separately
w*-cont}, we have
$$
b^{**}\widetilde{w}(\nu)=\lim_\alpha b^{**}w(f_\alpha) =
\lim_\alpha w(f_\alpha)b^{**} =\widetilde{w}(\nu)b^{**},
$$
and the claim follows.

Now fix some $f_1,\ldots,f_n\in C(K)$ and $b_1,\ldots,b_n\in E_u.$
For any $m\in \N$, there is a finite family $(t_1,\ldots, t_N)$ of
$K$ and a measurable partition $(B_1,\ldots, B_N)$ of $K$ such
that
$$
\Bignorm{f_k\, -\,\sum_{l=1}^N f_k(t_l)
\chi_{B_l}}_\infty\,\leq\,\frac{1}{m}\, ,\qquad k=1,\ldots,n.
$$
We set $f_{k}^{(m)} = \sum_{l=1}^N f_k(t_l) \chi_{B_l}$. Let
$\psi\colon \ell^\infty_N \to \BI(K)$ be the homomorphism of norm
1 mapping $e_l$ to $\chi_{B_l}$ for any $l$. According to
(\ref{seconddual}), we can consider the bounded homomorphism
$$
\widetilde{w}\circ \psi\colon \ell^\infty_N \longrightarrow
B(X\bi).
$$
Applying Proposition \ref{Prop Main thm for l infty n->B(X)} to
that homomorphism, together with the above claim and Lemma
\ref{Lem R(tau) = R(tau**)}, we obtain that
\begin{align*}
\Bignorm{\sum_k \widetilde{w}\bigl(f_k^{(m)}\bigr) b_k\bi} & =
\Bignorm{\sum_{k,l} f_k(t_l)\, \widetilde{w}\circ \psi(e_l)
b_k\bi}
\\
& \leq \|\widetilde{w}\circ \psi\|^{2} \, R\Bigr(\Bigl\{\sum_k
f_k(t_l) b_k\bi\, :\, 1\leq l\leq N\Bigr\}\Bigr)
\\
&\leq \|u\|^2 R\Bigr(\Bigl\{\sum_k f_k (t) b_k\bi\, :\, t\in
K\Bigr\}\Bigr)
\\ &\leq\norm{u}^{2}\,\Bignorm{\sum_k f_k\otimes b_k}_R.
\end{align*}
Since $\norm{f_k^{(m)}-f_k}_\infty \to 0$ for any $k$, we have
$$
\Bignorm{\sum_k \widetilde{w}\bigl(f_k^{(m)}\bigr) b_k\bi}
\longrightarrow \Bignorm{\sum_k w(f_k) b_k\bi} =\Bignorm{\sum_k
u(f_k) b_k},
$$
and the result follows at once.
\end{proof}

The following notion is implicit in several recent papers on
functional calculi (see in particular \cite{KaW,dPR}).

\begin{defi}\label{R-bounded map}
Let $Z$ be a Banach space and let $v\colon Z\to B(X)$ be a bounded
map. We set
$$
R(v)\,=\, R\bigl(\{v(z)\, :\, z\in Z,\ \norm{z}\leq 1\}\bigr),
$$
and we say that $v$ is $R$-bounded if $R(v)<\infty$.
\end{defi}

\begin{cor}\label{Cor Main Thm}
Let $u\colon C(K)\to B(X)$ be a bounded homomorphism and let
$v\colon Z\to B(X)$ be an $R$-bounded map. Assume further that
$u(f)v(z)=v(z)u(f)$ for any $f\in C(K)$ and any $z\in Z$. Then
there exists a (necessarily unique) bounded linear map
$$
u\cdotp v \colon C(K;Z) \longrightarrow B(X)
$$
such that $u\cdotp v(f\otimes z)=u(f)v(z)$ for any  $f\in C(K)$
and any $z\in Z$. Moreover we have
$$
\| u\cdotp v\|\leq \|u\|^2 R(v).
$$
\end{cor}

\begin{proof}
Consider any finite families $(f_k)_k$ in $C(K)$ and $(z_k)_k$ in
$Z$ and observe that
$$
\Bignorm{\sum_k f_k \otimes v(z_k)}_R =
R\Bigl(\Bigl\{v\Bigl(\sum_k f_k(t)z_k\Bigr)\, :\, t\in
K\Bigr\}\Bigr)\leq\,R(v)\,\Bignorm{\sum_k f_k \otimes z_k}_\infty.
$$
Then applying Theorem \ref{Thm Main} and the assumption that $v$
is valued in $E_u$, we obtain that
$$
\Bignorm{\sum_k u(f_k)v(z_k)}\leq \norm{u}^{2}
R(v)\,\Bignorm{\sum_k f_k \otimes z_k}_\infty,
$$
which proves the result.
\end{proof}

\begin{rem}\label{Rem De Pagter Ricker C(K) nuclear}
As a special case of Corollary \ref{Cor Main Thm}, we obtain the
following result due to De Pagter and Ricker (\cite[Prop.
2.27]{dPR}): Let $K_1,K_2$ be two compact sets, let
$$
u\colon C(K_1)\longrightarrow B(X)\qquad\hbox{and}\qquad v\colon
C(K_2)\longrightarrow  B(X)
$$
be two bounded homomorphisms which commute, i.e. $u(f)v(g) =
v(g)u(f)$ for all $f\in C(K_1)$ and $g\in C(K_2).$ Assume further
that $R(v)<\infty$. Then there exists a bounded homomorphism
$$
w\colon  C(K_1\times K_2)\longrightarrow B(X)
$$
such that $w_{|C(K_1)} = u$ and $w_{|C(K_2)} = v,$ where $C(K_j)$
is regarded as a subalgebra of $C(K_1\times K_2)$ in the natural
way.
\end{rem}

\medskip

\section{Uniformly bounded $\HI$-calculus}\label{Sec HI}
We briefly recall the basic notions on $\HI$-calculus for
sectorial operators. For more information, we refer  e.g. to
\cite{CDMY, KaW, KW, LM}.

For any $\theta\in (0,2\pi)$, we define
$$
\Sigma_\theta  = \{re^{i\phi}\, :\, r>0,\,|\phi|<\theta\}
$$
and
$$
\HI(\Sigma_\theta)  = \{f\colon\Sigma_\theta \to \C\,\vert\,
f\text{ is analytic and bounded}\}.
$$
This space is equipped with the norm $\|f\|_{\infty,\theta} =
\sup_{\lambda\in \Sigma_\theta} |f(\lambda)|$ and this is a Banach
algebra. We consider the auxiliary space
$$
\HI_0(\Sigma_\theta)= \bigl\{f \in \HI(\Sigma_\theta)\, :\,
\exists \: \epsilon,C>0\ \vert\ |f(\lambda)| \leq C
\min(|\lambda|^\epsilon,|\lambda|^{-\epsilon})\bigr\}.
$$
An closed linear operator $A\colon D(A)\subset X \to X$ is called
$\omega$-sectorial, for some $\omega\in(0,2\pi)$, if its domain
$D(A)$ is dense in $X$, its spectrum $\sigma(A)$ is contained in
$\overline{\Sigma_{\omega}}$, and for all $\theta
> \omega$ there is a constant $C_\theta >0$ such that
$$
\|\lambda (\lambda-A)^{-1}\| \leq C_\theta,\qquad \lambda \in
\C\setminus \overline{\Sigma_{\theta}}.
$$
In this case, we define
$$
\omega(A) = \inf\{\omega\, :\  A\text{ is
}\omega\text{-sectorial}\}.
$$
For any $\theta\in (\omega(A),\pi)$ and any $f\in
\HI_0(\Sigma_\theta)$, we define
\begin{equation}\label{CF}
f(A)\, =\, \frac{1}{2\pi i}\,\int_{\Gamma_\gamma} f(\lambda)
(\lambda-A)^{-1}\,d\lambda,
\end{equation}
where $\omega(A)<\gamma<\theta$ and $\Gamma_\gamma$ is the
boundary $\partial\Sigma_\gamma$ oriented counterclockwise. This
definition does not depend on $\gamma$ and the resulting mapping
$f\mapsto f(A)$ is an algebra homomorphism from
$\HI_0(\Sigma_\theta)$ into $B(X)$. We say that $A$ has a bounded
$\HI(\Sigma_\theta)$-calculus if the latter homomorphism is
bounded, that is, there exists a constant $C>0$ such that
$\|f(A)\|\leq C \|f\|_{\infty,\theta}$ for all $f\in
\HI_0(\Sigma_\theta).$

We  will now focus on sectorial operators $A$ such that
$\omega(A)=0$.

\begin{defi}\label{Uniform}
We say that a sectorial operator $A$ with $\omega(A) = 0$ has a
uniformly bounded $\HI$-calculus, if there exists a constant $C>0$
such that $\|f(A)\|\leq C\|f\|_{\infty,\theta}$ for all $\theta>0$
and $f\in \HI_0(\Sigma_\theta)$.
\end{defi}

We let
$$
C_\ell([0,\infty))=\bigl\{ f\colon [0,\infty)\to\C\, \vert \, f \
\hbox{is continuous and}\ \lim_{\infty}f\ \hbox{exists}\bigr\}.
$$
Then we equip this space with the sup norm
$$
\norm{g}_{\infty,0} =\sup\{\vert g(t)\vert\, :\, t>0\}.
$$
Thus $C_\ell([0,\infty))$ is a unital commutative $C^*$-algebra.
For any $\theta>0$, we can regard $\HI_0(\Sigma_\theta)$ as a
subalgebra of $C_\ell([0,\infty))$, by identifying any $f\in
\HI_0(\Sigma_\theta)$ with its restriction $f_{\vert [0,\infty)}$.

For any $\lambda\in\C\setminus [0,\infty)$, we let $R_\lambda\in
C_\ell([0,\infty))$ be defined by $R_\lambda(t)=(\lambda-t)^{-1}$.
Then we let $\mR$ be the unital algebra generated by the
$R_\lambda$'s. Equivalently, $\mR$ is the algebra of all rational
functions of nonpositive degree, whose poles lie outside the half
line $[0,\infty)$. We recall that for any $f\in
\HI_0(\Sigma_\theta)\cap\mR$, the definition of $f(A)$ given by
(\ref{CF}) coincides with the usual rational functional calculus.

\begin{lem}\label{Lem Technical HI}
Let $A$ be a sectorial operator on $X$ with $\omega(A) = 0.$ The
following assertions are equivalent.
\begin{enumerate}
\item [(a)] $A$ has a uniformly bounded $\HI$-calculus. \item
[(b)] There exists a (necessarily unique) bounded unital
homomorphism
$$
u\colon C_\ell([0,\infty)) \longrightarrow B(X)
$$
such that $u(R_\lambda) = (\lambda-A)^{-1}$ for any
$\lambda\in\C\setminus [0,\infty).$
\end{enumerate}
\end{lem}

\begin{proof}
Assume (a). We claim that for any $\theta>0$ and any $f\in
\HI_0(\Sigma_\theta)$, we have
$$
\norm{f(A)}\leq C\norm{f}_{\infty,0}.
$$
Indeed, if $0\neq f\in \HI_0(\Sigma_{\theta_0})$ for some
$\theta_0>0,$ then there exists some $t_0 > 0$ such that $f(t_0)
\neq 0.$ Now let $r<R$ such that $|f(z)| < |f(t_0)|$ for $|z|<r$
and $|z|>R.$ Choose for every $n\in \N$ a $t_n\in
\Sigma_{\theta_0/n}$ such that
$|f(t_n)|=\|f\|_{\infty,\theta_0/n}.$ Necessarily  $|t_n|\in
[r,R],$ and there exists a convergent subsequence $t_{n_k},$ whose
limit $t_\infty$ is real. Then
$$
\|f\|_{\infty,0}\geq |f(t_\infty)| \geq \liminf_{\theta\to 0}
\|f\|_{\infty,\theta} \geq C^{-1} \|f(A)\|.
$$
This readily implies that the rational functional calculus
$(\mR,\norm{\ }_{\infty,0})\to B(X)$ is bounded. By
Stone-Weierstrass, this extends continuously to
$C_\ell([0,\infty))$, which yields (b). The uniqueness property is
clear.

\smallskip
Assume (b). Then for any $\theta\in (0,\pi)$ and for any $f\in
\HI_0(\Sigma_{\theta})\cap\mR$, we have
$$
\norm{f(A)}\leq\norm{u}\norm{f}_{\infty,\theta}.
$$
By \cite[Prop. 2.10]{LM} and its proof, this implies that $A$ has
a bounded $\HI(\Sigma_\theta)$-calculus, with a boundedness
constant uniform in $\theta.$
\end{proof}

\begin{rem} An operator $A$ which admits a bounded
$\HI(\Sigma_\theta)$-calculus for all $\theta>0$ does not
necessarily have a uniformly bounded $\HI$-calculus. To get a
simple example, consider
$$
A = \left(
\begin{array}{cc} 1 & 1\\ 0 & 1 \\ \end{array}
\right)\,\colon\ell^{2}_{2}\longrightarrow \ell^{2}_{2}.
$$
Then $\sigma(A)=\{1\}$ and for any $\theta>0$ and any
$f\in\HI_{0}(\Sigma_\theta)$, we have
$$
f(A) = \left(
\begin{array}{cc} f(1) & f'(1)\\ 0 & f(1) \\ \end{array}
\right).
$$
Assume that $\theta<\frac{\pi}{2}$. Using Cauchy's Formula, it is
easy to see that $\vert
f'(1)\vert\leq(\sin(\theta))^{-1}\norm{f}_{\infty,\theta}$ for any
$f\in\HI_{0}(\Sigma_\theta)$. Thus $A$ admits a bounded
$\HI(\Sigma_\theta)$-calculus.

Now let $h$ be a fixed function in
$\HI_{0}(\Sigma_{\frac{\pi}{2}})$ such that $h(1)=1$, set
$g_s(\lambda) = \lambda^{is}$ for any $s>0$, and let $f_s=hg_s$.
Then $\norm{g_s}_{\infty,0}=1$, hence $\norm{f_s}_{\infty,0}\leq
\norm{h}_{\infty,0}$ for any $s>0$. Further $g'_s(\lambda) =
is\lambda^{is-1}$ and $f'_s=h'g_s + hg'_s$. Hence $f'_s(1) = h'(1)
+is$. Thus
$$
\norm{f_s(A)}\norm{f_s}_{\infty,0}^{-1}\geq \vert
f'_s(1)\vert\norm{h_s}_{\infty,0}^{-1}\longrightarrow \infty
$$
when $s\to\infty$. Hence $A$ does not have a uniformly bounded
$\HI$-calculus.

The above result can also be deduced from Proposition
\ref{Hilbert} below. In fact we will show in that proposition and
in Corollary \ref{Scalar} that operators with a uniformly bounded
$\HI$-calculus are `rare'.
\end{rem}

We now turn to the so-called generalized (or operator valued)
$\HI$-calculus. Throughout we let $A$ be a sectorial operator. We
let $E_A\subset B(X)$ denote the commutant of $A$, defined as the
subalgebra of all bounded operators $T\colon X\to X$ such that
$T(\lambda-A)^{-1}= (\lambda-A)^{-1}T$ for any $\lambda$ belonging
to the resolvent set of $A$. We let $\HI_0(\Sigma_\theta; B(X))$
be the algebra of all bounded analytic functions $F
\colon\Sigma_\theta\to B(X)$ for which there exist $\epsilon,C>0$
such that $\norm{F(\lambda)} \leq C
\min(|\lambda|^\epsilon,|\lambda|^{-\epsilon})$ for any
$\lambda\in\Sigma_\theta$. Also, we let $\HI_0(\Sigma_\theta;
E_A)$ denote the space of all $E_A$-valued functions belonging to
$\HI_0(\Sigma_\theta; B(X))$. The generalized $\HI$-calculus of
$A$ is an extension of (\ref{CF}) to this class of functions.
Namely for any $F\in \HI_0(\Sigma_\theta; E_A)$, we set
$$
F(A)\, =\, \frac{1}{2\pi i}\,\int_{\Gamma_\gamma} F(\lambda)
(\lambda-A)^{-1}\,d\lambda\,,
$$
where $\gamma\in(\omega(A),\pi)$. Again, this definition does not
depend on $\gamma$ and the mapping $F\mapsto F(A)$ is an algebra
homomorphism. The following fundamental result is due to Kalton
and Weis.

\begin{thm}\label{Thm Kalton Weis Operator Valued HI}
(\cite[Thm. 4.4]{KaW},\,\cite[Thm. 12.7]{KW}.) Let $\omega_0\geq
\omega(A)$ and assume that $A$ has a bounded
$\HI(\Sigma_\theta)$-calculus for any $\theta>\omega_0$. Then for
any $\theta>\omega_0$, there exists a constant $C_\theta>0$ such
that for any $F\in \HI_0(\Sigma_\theta;E_A)$,
\begin{equation}\label{KW}
\norm{F(A)} \leq C_\theta R\bigl(\{F(z)\, :\
z\in\Sigma_\theta\}\bigr).
\end{equation}
\end{thm}

Our aim is to prove a version of this result in the case when $A$
has a uniformly bounded $\HI$-calculus. We will obtain in Theorem
\ref{Thm HI Main} that in this case, the constant $C_\theta$ in
(\ref{KW}) can be taken independent of $\theta$.

The algebra $C_\ell([0,\infty))$ is a $C(K)$-space and we will
apply the results of Section 2 to the bounded homomophism $u$
appearing in Lemma \ref{Lem Technical HI}.  We recall Remark
\ref{Embed}.

\begin{lem}\label{Automatic R} Let $J\colon C_\ell( [0,\infty))
\mathop{\otimes}\limits^R B(X)\to
C_\ell([0,\infty);B(X))$ be the canonical embedding. Let
$\theta\in(0,\pi)$, let $F\in \HI_0(\Sigma_\theta; B(X))$ and let
$\gamma\in(0,\theta)$.
\begin{enumerate}
\item [(1)] The integral
\begin{equation}\label{Integral}
\varphi_F\, =\, \frac{1}{2\pi i}\,\int_{\Gamma_\gamma}
R_\lambda\otimes F(\lambda)\, d\lambda
\end{equation}
is absolutely convergent in
$C_\ell([0,\infty))\mathop{\otimes}\limits^R B(X)$, and
$J(\varphi_F)$ is equal to the restriction of $F$ to $[0,\infty)$.
\item [(2)] The set $\{F(t)\, :\, t>0\}$ is $R$-bounded.
\end{enumerate}
\end{lem}

\begin{proof}
Part (2) readily follows from part (1) and Remark \ref{Embed}. To
prove (1), observe that for any
$\lambda\in\partial\Sigma_{\gamma}$, we have
$$
\norm{R_\lambda\otimes F(\lambda)}_R\leq
2\norm{R_\lambda}_{\infty,0}\norm{F(\lambda)}
\leq\,\frac{2}{\sin(\gamma)\vert\lambda\vert}\, \norm{F(\lambda)}
$$
by (\ref{R}). Thus for appropriate constants $\epsilon, C>0$, we
have
$$
\norm{R_\lambda\otimes F(\lambda)}_R\leq \frac{2C}{\sin(\gamma)}\,
\min(\vert\lambda\vert^{\epsilon
-1},\vert\lambda\vert^{-\epsilon-1}).
$$
This shows that the integral defining $\varphi_F$ is absolutely
convergent. Next, for any $t>0$, we have
$$
[J(\varphi_F)](t)\, =\, \frac{1}{2\pi i}\,\int_{\Gamma_\gamma}
\bigl(R_\lambda \otimes F(\lambda)\bigr)(t)\, d\lambda\, = \,
\frac{1}{2\pi i}\,\int_{\Gamma_\gamma} \frac{F(\lambda)}{\lambda
-t}\, d\lambda\, =\, F(t)
$$
by Cauchy's Theorem.
\end{proof}

\begin{thm}\label{Thm HI Main} Let $A$ be a sectorial operator
with $\omega(A)=0$ and assume that $A$ has a uniformly bounded
$\HI$-calculus. Then there exists a constant $C>0$ such that for
any $\theta>0$ and any $F\in \HI_0(\Sigma_\theta;E_A)$,
$$
\|F(A)\|\leq C R\bigl(\{F(t)\, :\, t>0\}\bigr).
$$
\end{thm}

\begin{proof}
Let $u\colon C_\ell([0,\infty))\to B(X)$ be the representation
given by Lemma \ref{Lem Technical HI}. It is plain that $E_u=E_A$.
Then we let
$$
\widehat{u}\colon C_\ell([0,\infty)) \mathop{\otimes}\limits^R
E_A\longrightarrow B(X)
$$
be the associated bounded map provided by Theorem \ref{Thm Main}.

Let $F\in \HI_0(\Sigma_\theta; E_A)$ for some $\theta>0$, and let
$\varphi_F\in C_\ell([0,\infty))\mathop{\otimes}\limits^R E_A$ be
defined by (\ref{Integral}). We claim that
$$
F(A)=\widehat{u}(\varphi_F).
$$
Indeed for any $\lambda\in\partial\Sigma_\gamma$, we have
$u(R_\lambda)=(\lambda -A)^{-1}$, hence
$\widehat{u}\bigl(R_\lambda\otimes F(\lambda)\bigr) =
(\lambda-A)^{-1}F(\lambda)$. Thus according to the definition of
$\varphi_F$ and the continuity of $\widehat{u}$, we have
$$
\widehat{u}(\varphi_F)=\, \frac{1}{2\pi i}\,\int_{\Gamma_\gamma}
\widehat{u}\bigl(R_\lambda\otimes F(\lambda)\bigr)\, d\lambda\ =
\frac{1}{2\pi i}\,\int_{\Gamma_\gamma} (\lambda-A)^{-1}
F(\lambda)\, d\lambda\ = F(A).
$$
Consequently,
$$
\norm{F(A)}\leq\norm{\widehat{u}}\norm{\varphi_F}_R\leq
\norm{u}^2\norm{\varphi_F}_R.
$$
It follows from Lemma \ref{Automatic R} and Remark \ref{Embed}
that $\norm{\varphi_F}_R =R\bigl(\{F(t)\, :\, t>0\}\bigr)$, and
the result follows at once.
\end{proof}

In the rest of this section we will further investigate operators
with a uniformly bounded $\HI$-calculus. We start with the case
when $X$ is a Hilbert space.

\begin{prop}\label{Hilbert}
Let $H$ be a Hilbert space and let $A$ be a sectorial operator on
$H$, with $\omega(A)=0$. Then $A$ admits a uniformly bounded
$\HI$-calculus if and only if there exists an isomorphism $S\colon
H\to H$ such that $S^{-1}AS$ is selfadjoint.
\end{prop}

\begin{proof}
Assume that $A$ admits a uniformly bounded $\HI$-calculus and let
$u\colon C_\ell([0,\infty))\to B(H)$ be the associated
representation. According to \cite[Thm. 9.1 and Thm. 9.7]{Pa},
there exists an isomorphism $S\colon H\to H$ such that the unital
homomorphism $u_S\colon C_\ell([0,\infty))\to B(H)$ defined by
$u_S(f) = S^{-1}u(f)S$ satisfies $\norm{u_S}\leq 1$. We let
$B=S^{-1}AS$. For any $s\in\R^{*}$, we have
$\norm{R_{is}}_{\infty,0}=\vert s\vert$ and $u_S(R_{is}) =
S^{-1}(is-A)^{-1}S = (is-B)^{-1}$. Hence
$$
\norm{(is -B)^{-1}}\leq \vert s\vert,\qquad s\in\R^{*}.
$$
By the Hille-Yosida Theorem, this implies that $iB$ and $-iB$ both
generate contractive $c_0$-semigroups on $H$. Thus $iB$ generates
a unitary $c_0$-group. By Stone's Theorem, this implies that $B$
is selfadjoint.

The converse implication is clear.
\end{proof}

In the non Hilbertian setting, we will first show that operators
with a uniformly bounded $\HI$-calculus satisfy a spectral mapping
theorem with respect to continuous functions defined on the
one-point compactification of $\sigma(A)$. Then we will discuss
the connections with spectral measures and scalar-type operators.
We mainly refer to \cite[Chap. 5-7]{Dow} for this topic.

For any compact set $K$ and any closed subset $F\subset K$, we let
$$
I_F = \{f\in C(K)\, :\, f_{|F} = 0\}.
$$
We recall that the restriction map $f\mapsto f_{\vert F}$ induces
a $*$-isomorphism $C(K)/I_F\to C(F)$.

\begin{lem}\label{Lem Hom Spectrum Reduction}
Let $K\subset \C$ be a compact set and let $u\colon C(K) \to B(X)$
be a representation. Let $\kappa \in C(K)$ be the function defined
by $\kappa(z) = z$ and put $T = u(\kappa).$
\begin{enumerate}
\item [(1)] Then $\sigma(T) \subset K$ and $u$ vanishes on
$I_{\sigma(T)}$.
\end{enumerate}
Let $v\colon C(\sigma(T))\simeq C(K)/I_{\sigma(T)}\longrightarrow
B(X)$ be the representation induced by $u$.
\begin{enumerate}
\item [(2)] For any $f\in C(\sigma(T)),$ we have $\sigma(v(f)) =
f(\sigma(T)).$ \item [(3)] $v$ is an isomorphism onto its range.
\end{enumerate}
\end{lem}

\begin{proof} The inclusion $\sigma(T)\subset K$ is clear.
Indeed, for any $\lambda\notin K$, $(\lambda-T)^{-1}$ is equal to
$u((\lambda -\,\cdotp)^{-1})$. We will now show that $u$ vanishes
on $I_{\sigma(T)}$.

Let $w\colon C(K)\to B(X^{*})$ be defined by $w(f)=[u(f)]^{*}$,
and let $\widetilde{w} \colon C(K)\bi \to B(X^*)$ be its
$w^*$-extension. Since $w$ is valued in $w^*B(X^*)\simeq B(X)$,
this is a representation (see Lemma \ref{Lem Extension C(K)->B(X)
to C(K)**->B(X)}). Let $\Delta_K$ be the set of all Borel subsets
of $K$. It is easy to check that the mapping
$$
P \colon\Delta_K \longrightarrow B(X^*),\qquad P(B) =
\widetilde{w}(\chi_B),
$$
is a spectral measure of class $(\Delta_K ,X)$ in the sense of
\cite[p. 119]{Dow}. According to \cite[Prop. 5.8]{Dow}, the
operator $T^*$ is prespectral of class $X$ (in the sense of
\cite[Def. 5.5]{Dow}) and the above mapping $P$ is its resolution
of the identity. Applying \cite[Lem. 5.6]{Dow} and the equality
$\sigma(T^*)=\sigma(T)$, we obtain that
$\widetilde{w}(\chi_{\sigma(T)}) = P(\sigma(T)) = I_{X^{*}}.$
Therefore  for any $f\in I_{\sigma(T)}$, we have
$$
u(f)^{*} = \widetilde{w}\bigl(f(1-\chi_{\sigma(T)})\bigr) =
\widetilde{w}(f)\widetilde{w}(1-\chi_{\sigma(T)}) = 0.
$$
Hence $u$ vanishes on $I_{\sigma(T)}$.

The proofs of (2) and (3) now follow from \cite[Prop. 5.9]{Dow}
and the above proof.
\end{proof}

In the sequel we consider a sectorial operator $A$ with
$\omega(A)=0$. This assumption implies that
$\sigma(A)\subset[0,\infty)$. By $C_\ell(\sigma(A))$, we denote
either the space $C(\sigma(A))$ if $A$ is bounded, or the space
$\{f\colon \sigma(A)\to\C\, \vert \,f\ \hbox{is continuous and}\
\lim_\infty f \ \hbox{exists}\,\}$ if $A$ is unbounded. In this
case, $C_\ell(\sigma(A))$ coincides with the space of continuous
functions on the one-point compactification of $\sigma(A)$. The
following strengthens Lemma \ref{Lem Technical HI}.

\begin{prop}\label{Thm Description uniform HI}
Let $A$ be a sectorial operator on $X$ with $\omega(A) = 0.$ The
following assertions are equivalent.
\begin{enumerate}
\item[(1)] $A$ has a uniformly bounded $\HI$-calculus. \item[(2)]
There exists a (necessarily unique) bounded unital homomorphism
$$
\Psi\colon C_\ell(\sigma(A)) \longrightarrow  B(X)
$$
such that $\Psi\bigl((\lambda-\,\cdotp)^{-1}\bigr) =
(\lambda-A)^{-1}$ for any $\lambda \in \C\setminus\sigma(A)$.
\end{enumerate}
In this case, $\Psi$ is an isomorphism onto its range and for any
$f\in C_\ell(\sigma(A))$, we have
\begin{equation}\label{Sp}
\sigma(\Psi(f)) = f(\sigma(A))\cup f_\infty,
\end{equation}
where $f_\infty = \emptyset$ if $A$ is bounded and $f_\infty =
\{\lim_\infty f\}$ if $A$ is unbounded.
\end{prop}

\begin{proof}
Assume (1) and let $u\colon C_\ell([0,\infty))\to B(X)$ be given
by Lemma \ref{Lem Technical HI}. We introduce the particular
function $\phi\in C_\ell([0,\infty))$ defined by $\phi(t) =
(1+t)^{-1}.$ Consider the $*$-isomorphism
$$
\tau\colon C([0,1]) \longrightarrow C_\ell([0,\infty)),\qquad
\tau(g) = g\circ \phi,
$$
and set $T=(1+A)^{-1}$. If we let $\kappa(z)=z$ as in Lemma
\ref{Lem Hom Spectrum Reduction}, we have
$(u\circ\tau)(\kappa)=T$. Let $v\colon  C(\sigma(T))\to B(X)$ be
the resulting factorisation of $u \circ \tau$. The spectral
mapping theorem gives $\sigma(A) = \phi^{-1}(\sigma(T)\setminus
\{0\})$ and $0\in \sigma(T)$ if and only if $A$ is unbounded. Thus
the mapping
$$
\tau_A \colon C(\sigma(T)) \longrightarrow C_\ell (\sigma(A))
$$
defined by $\tau_A(g)= g\circ \phi$ also is a $*$-isomorphism. Put
$\Psi = v \circ \tau_{A}^{-1} \colon C_\ell(\sigma(A))\to B(X)$.
This is a unital bounded homommorphism. Note that $\phi^{-1}(z) =
\frac{1-z}{z}$ for any $z\in(0,1]$. Then  for any $\lambda\in
\C\setminus\sigma(A)$,
\begin{align*}
\Psi\bigl((\lambda-\,\cdotp)^{-1}\bigr) & =
v\bigl((\lambda-\,\cdotp)^{-1}\circ\phi^{-1}) = v\biggl(z\mapsto
\Bigl(\lambda
-\,\frac{1-z}{z}\Bigr)^{-1}\biggr)\\
& = v\Bigl(z\mapsto \frac{z}{(\lambda+1)z -1}\Bigr)\\ & =
T\bigl((\lambda +1)T -1\bigr)^{-1} = (\lambda -A)^{-1}.
\end{align*}
Hence $\Psi$ satisfies (2). Its uniqueness follows from Lemma
\ref{Lem Technical HI}. The fact that $\Psi$ is an isomorphism
onto its range, and the spectral property (\ref{Sp}) follow from
the above construction and Lemma \ref{Lem Hom Spectrum Reduction}.
Finally the implication `(2)$\,\Rightarrow\,$(1)' also follows
from Lemma \ref{Lem Technical HI}.
\end{proof}

\begin{rem}
Let $A$ be a sectorial operator with a uniformly bounded
$\HI$-calculus, and let $T=(1+A)^{-1}$. It follows from Lemma
\ref{Lem Hom Spectrum Reduction} and the proof of Proposition
\ref{Thm Description uniform HI} that there exists a
representation
$$
v\colon C(\sigma(T))\longrightarrow B(X)
$$
satisfying $v(\kappa)=T$ (where $\kappa(z)=z$), such that
$\sigma(v(f))=f(\sigma(T))$ for any $f\in C(\sigma(T))$ and $v$ is
an isomorphism onto its range. Also, it follows from the proof of
Lemma \ref{Lem Hom Spectrum Reduction} that $T^*$ is a scalar-type
operator of class $X$, in the sense of \cite[Def. 5.14]{Dow}.

Next according to \cite[Thm. 6.24]{Dow}, the operator $T$ (and
hence $A$) is a scalar-type spectral operator if and only if for
any $x\in X$, the mapping $C(\sigma(T))\to X$ taking $f$ to
$v(f)x$ for any $f\in C(\sigma(T))$ is weakly compact.
\end{rem}

\begin{cor}\label{Scalar} Let $A$ be a sectorial operator on $X$, with
$\omega(A)=0$, and assume that $X$ does not contain a copy of
$c_0$. Then $A$ admits a uniformly bounded $\HI$-calculus if and
only if it is a scalar-type spectral operator.
\end{cor}

\begin{proof}
The `only if' part follows from the previous remark. Indeed if $X$
does not contain a copy of $c_0$, then any bounded map $C(K)\to X$
is weakly compact \cite[VI, Thm. 15]{DU}. (See also \cite{R} and
\cite{dPR} for related approaches.) The `if' part  follows from
\cite[Prop. 2.7]{FW} and its proof.
\end{proof}

\begin{rem}
Scalar-type spectral operator on Hilbert space coincide with
operators similar to a normal one (see \cite[Chap. 7]{Dow}). Thus
when $X=H$ is a Hilbert space, the above corollary reduces to
Proposition \ref{Hilbert}.
\end{rem}

\medskip

\section{Matricial $R$-boundedness}\label{Sec Matr R-bound}

For any integer $n\geq 1$ and any vector space $E$, we will denote
by $M_n(E)$ the space of all $n\times n$ matrices with entries in
$E$. We will be mostly concerned with the cases $E=C(K)$ or
$E=B(X)$. As mentioned in the introduction, we identify
$M_n(C(K))$ with the space $C(K;M_n)$ in the usual way. We now
introduce a specific norm on $M_n(B(X))$. Namely for any
$[T_{ij}]\in M_n(B(X))$, we set
$$
\bignorm{[T_{ij}]}_R\, =\,
\sup\Bigl\{\Bignorm{\sum_{i,j=1}^{n}\epsilon_i\otimes
T_{ij}(x_j)}_{{\rm Rad}(X)}\, :\, x_1,\ldots, x_n\in X, \
\Bignorm{\sum_{j=1}^{n}\epsilon_j\otimes  x_j}_{{\rm Rad}(X)}\leq
1\,\Bigr\}.
$$
Clearly $\norm{\ }_R$ is a norm on $M_n(B(X))$. Moreover if we
consider any element of $M_n(B(X))$ as an operator on
$\ell^2_n\otimes X$ in the natural way, and if we equip the latter
tensor product with the norm of ${\rm Rad}_{n}(X)$, we obtain an
isometric identification
\begin{equation}\label{Ident}
\bigl(M_n(B(X)),\norm{\ }_R\bigr)\,=\, B\bigl({\rm
Rad}_{n}(X)\bigr).
\end{equation}

\begin{defi}\label{Def matricially R-bounded}
Let $u\colon C(K)\to B(X)$ be a bounded linear mapping. We say
that $u$ is matricially $R$-bounded if there is a constant $C\geq
0$ such that for any $n\geq 1$, and for any $[f_{ij}]\in
M_n(C(K))$, we have
\begin{equation}\label{Mat}
\bignorm{[u(f_{ij})]}_R \,\leq\, C \norm{[f_{ij}]}_{C(K;M_n)}.
\end{equation}
\end{defi}

\begin{rem}\label{4Remark}
The above definition obviously extends to any bounded map $E\to
B(X)$ defined on an operator space $E$, or more generally on any
matricially normed space (see \cite{ER1,ER2}). The basic
observations below apply to this general case as well.

\smallskip
(1) In the case when $X=H$ is a Hilbert space, we have
$$
\Bignorm{\sum_{j=1}^{n}\epsilon_j\otimes  x_j}_{{\rm
Rad}(H)}\,=\,\Bigl(\sum_{j=1}^{n}\norm{x_j}^2\Bigr)^{\frac{1}{2}}
$$
for any $x_1,\ldots, x_n\in H$. Consequently, writing that a
mapping $u\colon C(K)\to B(H)$ is matricially $R$-bounded is
equivalent to writing that $u$ is completely bounded (see e.g.
\cite{Pa}). See Section 5 for the case when $X$ is an $L^p$-space.

\smallskip (2) The notation $\norm{\ }_R$ introduced above is
consistent with the one considered so far in Section 2. Indeed let
$b_1,\ldots,b_n$ in $B(X)$. Then the diagonal matrix ${\rm
Diag}\{b_1,\ldots,b_n\}\in M_n(B(X))$ and the tensor element
$\sum_{k=1}^{n} e_k\otimes b_k\,\in\ell^{\infty}_n\otimes B(X)$
satisfy
$$
\bignorm{{\rm
Diag}\{b_1,\ldots,b_n\}}_R\,=\,R(\{b_1,\ldots,b_n\})\,=\,\Bignorm{\sum_{k=1}^{n}
e_k\otimes b_k}_R.
$$

\smallskip (3) If $u\colon C(K)\to B(X)$ is
matricially $R$-bounded (with the estimate (\ref{Mat})), then $u$
is $R$-bounded and $R(u)\leq C$. Indeed, consider $f_1,\ldots,
f_n$ in the unit ball of $C(K)$. Then we have $\norm{{\rm
Diag}\{f_1,\ldots,f_n\}}_{C(K;M_n)}\leq 1$. Hence for any
$x_1,\ldots,x_n$ in $X$,
\begin{align*}
\Bignorm{\sum_k\epsilon_k\otimes u(f_k)x_k}_{{\rm Rad}(X)}\, &
\leq\,\bignorm{{\rm Diag}\{u(f_1),\ldots,u(f_n)\}}_R\,
\Bignorm{\sum_k\epsilon_k\otimes  x_k}_{{\rm Rad}(X)}\\
& \leq\, C\Bignorm{\sum_k\epsilon_k\otimes  x_k}_{{\rm Rad}(X)}.
\end{align*}
\end{rem}

Let $(g_k)_{k\geq 1}$ be a sequence of complex valued,
independent, standard Gaussian random variables on some
probability space $\Omega_G$. For any $x_1,\ldots, x_n$ in $X$ let
$$
\Bignorm{\sum_k g_k\otimes
x_k}_{G(X)}\,=\,\Bigl(\int_{\Omega_G}\Bignorm{\sum_k g_k(\lambda)
\, x_k}^{2}_{X}\, d\lambda\,\Bigr)^{\frac{1}{2}}.
$$
It is well-known that for any scalar valued matrix $a=[a_{ij}]\in
M_n$, we have
\begin{equation}\label{Invariant}
\Bignorm{\sum_{i,j=1}^{n} a_{ij}\, g_i \otimes
x_j}_{G(X)}\leq\,\norm{a}_{M_n}\,\Bignorm{\sum_{j=1}^{n}
g_j\otimes x_j}_{G(X)},
\end{equation}
see e.g. \cite[Cor. 12.17]{DJT}. For any $n\geq 1$, introduce
$$
\sigma_{n,X}\colon
\begin{cases}
\quad\  M_n & \longrightarrow B\bigl({\rm Rad}_{n}(X)\bigr)\\
a=[a_{ij}]  & \longmapsto [a_{ij} I_X].
\end{cases}
$$
If $X$ has finite cotype, then we have a uniform equivalence
\begin{equation}\label{cotype}
\Bignorm{\sum_k \epsilon_k\otimes x_k}_{{\rm Rad}(X)} \asymp
\Bignorm{\sum_k g_k\otimes x_k}_{G(X)}
\end{equation}
between Rademacher and Gaussian averages on $X$ (see e.g.
\cite[Thm. 12.27]{DJT}). Combining with (\ref{Invariant}), this
implies that
$$
\sup_{n\geq 1}\norm{\sigma_{n,X}}\,<\infty.
$$

Following \cite{Pis} we say that $X$ has property $(\alpha)$ if
there is a constant $C\geq 1$ such that for any finite family
$(x_{ij})$ in $X$ and any finite family $(t_{ij})$ of complex
numbers,
\begin{equation}\label{alpha}
\Bignorm{\sum_{i,j} \epsilon_i\otimes \epsilon_j \otimes
t_{ij}x_{ij}}_{{\rm Rad}({\rm Rad}(X))}\,\leq\, C\,
\sup_{i,j}\vert t_{ij}\vert\,\Bignorm{\sum_{i,j} \epsilon_i\otimes
\epsilon_j \otimes  x_{ij}}_{{\rm Rad}({\rm Rad}(X))}.
\end{equation}
Equivalently, $X$ has property $(\alpha)$ if and only if we have a
uniform equivalence
$$
\Bignorm{\sum_{i,j} \epsilon_i\otimes \epsilon_j \otimes
x_{ij}}_{{\rm Rad}({\rm Rad}(X))}\,\asymp\, \Bignorm{\sum_{i,j}
\epsilon_{ij} \otimes x_{ij}}_{{\rm Rad}(X)},
$$
where $(\epsilon_{ij})_{i,j\geq 1}$ be a doubly indexed family of
independent Rademacher variables.

The following is a characterization of property $(\alpha)$ in
terms of the $R$-boundedness of $\sigma_{n,X}$.

\begin{lem}\label{Lem M_n matrically R-bounded}
A Banach space $X$ has property $(\alpha)$ if and only if
$$
\sup_{n\geq 1} R(\sigma_{n,X})\,<\infty.
$$
\end{lem}

\begin{proof}
Assume that $X$ has property $(\alpha).$ This implies that $X$ has
finite cotype, hence $X$ satisfies the equivalence property
(\ref{cotype}). Let $a(1),\ldots, a(N)$ in $M_n$ and let
$z_1,\ldots, z_N$ in ${\rm Rad}_{n}(X)$. Let $x_{jk}$ in $X$ such
that $z_k=\sum_j\epsilon_j\otimes x_{jk}\,$ for any $k$. We
consider a doubly indexed family $(\epsilon_{ik})_{i,k\geq 1}$ as
above, as well as a doubly indexed family $(g_{ik})_{i,k\geq 1}$
of independent standard Gaussian variables. Then
\begin{equation}\label{def snX}
\sum_k\epsilon_k\otimes \sigma_{n,X}\bigl(a(k)\bigr)z_k\, =\,
\sum_{k,i,j} \epsilon_k\otimes \epsilon_i \otimes a(k)_{ij}
x_{jk}.
\end{equation}
Hence using the properties reviewed above, we have
\begin{align*}
\Bignorm{\sum_k\epsilon_k\otimes
\sigma_{n,X}\bigl(a(k)\bigr)z_k}_{{\rm Rad}({\rm Rad}(X))}\, &
\asymp\, \Bignorm{\sum_{k,i,j} \epsilon_{ik}\otimes a(k)_{ij}
x_{jk}}_{{\rm Rad}(X)}\\
& \asymp\, \Bignorm{\sum_{k,i,j} g_{ik}\otimes a(k)_{ij}
x_{jk}}_{G(X)}\\
& \lesssim\, \left\|\left(
\begin{array}{ccc}
 a(1) & 0 \ldots & 0 \\
 0 & \ddots & 0 \\
 0 & \ldots 0 & a(N) \\
\end{array}
\right)\right\|_{M_{Nn}}\,\Bignorm{\sum_{k,j} g_{jk}\otimes
x_{jk}}_{G(X)}\\
& \lesssim\,\max_k\norm{a(k)}_{M_n}\, \Bignorm{\sum_{k,j}
\epsilon_{jk}\otimes x_{jk}}_{{\rm Rad}(X)}\\
& \lesssim\,\max_k\norm{a(k)}_{M_n}\, \Bignorm{\sum_{k,j}
\epsilon_{k}\otimes \epsilon_j \otimes x_{jk}}_{{\rm Rad}({\rm Rad}(X))}\\
& \lesssim\,\max_k\norm{a(k)}_{M_n}\, \Bignorm{\sum_{k}
\epsilon_{k}\otimes z_k}_{{\rm Rad}({\rm Rad}(X))}.
\end{align*}
This shows that the $\sigma_{n,X}$'s are uniformly $R$-bounded.

Conversely, assume that for some constant $C\geq 1$, we have
$R(\sigma_{n,X})\leq C$ for any $n\geq 1$. Let $(t_{jk})_{j,k}\in
\C^{n^{2}}$ with $|t_{jk}|\leq 1$ and for any $k=1,\ldots,n$, let
$a(k)\in M_n$ be the diagonal matrix with entries
$t_{1k},\ldots,t_{nk}$ on the diagonal. Then $\norm{a(k)}\leq 1$
for any $k$. Hence applying (\ref{def snX}), we obtain that for
any $(x_{jk})_{j,k}$ in $X^{n^{2}}$,
\begin{align*}
\Bignorm{\sum_{j,k}\epsilon_k \otimes\epsilon_j\otimes
t_{jk}x_{jk}}_{\Rad(\Rad(X))}\, & \leq\,
R\bigl(\{a(1),\ldots,a(n)\}\bigr)\, \Bignorm{\sum_{j,k}\epsilon_k
\otimes\epsilon_j\otimes x_{jk}}_{\Rad(\Rad(X))}\\
& \leq\, C\, \Bignorm{\sum_{j,k}\epsilon_k
\otimes\epsilon_j\otimes x_{jk}}_{\Rad(\Rad(X))}.
\end{align*}
This means that $X$ has property $(\alpha).$
\end{proof}

\begin{prop}\label{Prop u matricially R-bounded}
Assume that $X$ has property $(\alpha).$ Then any bounded
homomorphism $u\colon C(K) \to B(X)$ is matricially $R$-bounded.
\end{prop}

\begin{proof}
Let $u\colon C(K) \to B(X)$ be a bounded homomorphism and let
$$
w \colon \begin{cases}
C(K) & \longrightarrow B\bigl(\Rad_{n}(X)\bigr) \\
f    & \longmapsto I_{\Rad_{n}}\otimes u(f).
\end{cases}
$$
Clearly $w$ is also a bounded homomorphism, with $\|w\|=\|u\|.$
Recall the identification (\ref{Ident}) and note that $w(f) = {\rm
Diag}\{u(f),\ldots,u(f)\}$ for any $f\in C(K)$. Then for any
$a=[a_{ij}]\in M_n$, we have
$$
w(f)\sigma_{n,X}(a) = [a_{ij}u(f)] = \sigma_{n,X}(a) w(f).
$$
By Corollary \ref{Cor Main Thm} and Lemma \ref{Lem M_n matrically
R-bounded}, the resulting mapping $w\cdotp \sigma_{n,X}$ satisfies
$$
\bignorm{w\cdotp \sigma_{n,X}\colon C(K; M_n) \longrightarrow
B(\Rad_{n}(X))} \leq C \|u\|^2
$$
where $C$ does not depend on $n.$ Let $E_{ij}$ denote the
canonical matrix units of $M_n$, for $i,j=1,\ldots,n$. Consider
$[f_{ij}]\in C(K;M_n)\simeq M_n(C(K))$ and write this matrix as
$\sum_{i,j} E_{ij}\otimes f_{ij}\,$. Then
$$
w\cdotp \sigma_{n,X}\bigl([f_{ij}]\bigr) =
\sum_{i,j=1}^{n}w(f_{ij})\sigma_{n,X}(E_{ij})\,= \sum_{i,j=1}^{n}
u(f_{ij})\otimes E_{ij}\,=[u(f_{ij})].
$$
Hence $\bignorm{[u(f_{ij})]}_R\leq
C\norm{u}^2\norm{[f_{ij}]}_{C(K;M_n)}$, which proves that $u$ is
matricially $R$-bounded.
\end{proof}

Suppose that $X=H$ is a Hibert space, and recall Remark
\ref{4Remark} (1). Then in that case, the above proposition
reduces to the fact that any bounded homomorphism $C(K)\to B(H)$
is completely bounded.

We also observe that applying the above proposition together with
Remark \ref{4Remark} (3), we obtain the following corollary
originally due to De Pagter and Ricker \cite[Cor. 2.19]{dPR}.
Indeed, Proposition \ref{Prop u matricially R-bounded} should be
regarded as a strengthening of their result.

\begin{cor}\label{Cor bounded (alpha) R-bounded}$\,$
Assume that $X$ has property $(\alpha)$. Then any bounded
homomorphism $u\colon C(K) \to B(X)$ is $R$-bounded.
\end{cor}

\begin{rem} The above corollary is nearly optimal.  Indeed we claim that if
$X$ does not have property $(\alpha)$ and if $K$ is any infinite
compact set, then there exists a unital bounded homomorphism
$$
u\colon C(K)\longrightarrow B\bigl({\rm Rad}(X)\bigr)
$$
which is not $R$-bounded.

To prove this, let $(z_n)_{n\geq 1}$ be an infinite sequence of
distinct points in $K$ and let $u$ be defined by
$$
u(f)\Bigl(\sum_{k\geq 1}\epsilon_k\otimes
x_k\Bigr)\,=\,\sum_{k\geq 1} f(z_k)\,\epsilon_k\otimes x_k.
$$
According to (\ref{CP}), this is a bounded unital homomorphism
satisfying $\norm{u}\leq 2$. Assume now that $u$ is $R$-bounded.
Let $n\geq 1$ be an integer and consider families $(t_{ij})_{i,j}$
in $\C^{n^{2}}$ and $(x_{ij})_{i,j}$ in $X^{n^{2}}$. For any
$i=1,\ldots, n$, there exists $f_i\in C(K)$ such that
$\norm{f_i}=\sup_j\vert t_{ij}\vert$ and $f_i(z_j) = t_{ij}$ for
any $j=1,\ldots, n$. Then
$$
\sum_{i}\epsilon_i\otimes u(f_i)\Bigl(\sum_j\epsilon_j\otimes
x_{ij}\Bigr)\,=\,\sum_{i,j}
t_{ij}\,\epsilon_i\otimes\epsilon_j\otimes x_{ij},
$$
hence
\begin{align*}
\Bignorm{\sum_{i,j} t_{ij}\,\epsilon_i\otimes\epsilon_j\otimes
x_{ij}}_{{\rm Rad}({\rm Rad}(X))}\, &\leq
R(u)\sup_i\norm{f_i}\,\Bignorm{\sum_{i,j}
\epsilon_i\otimes\epsilon_j\otimes x_{ij}}_{{\rm Rad}({\rm
Rad}(X))}\\ & \leq R(u)\sup_{i,j}\vert t_{ij}\vert\,
\Bignorm{\sum_{i,j} \epsilon_i\otimes\epsilon_j\otimes
x_{ij}}_{{\rm Rad}({\rm Rad}(X))}.
\end{align*}
This shows (\ref{alpha}).
\end{rem}

\medskip

\section{Application to $L^p$-spaces and unconditional bases}\label{Sec Bases}

Let $X$ be a Banach lattice with finite cotype. A classical
theorem of Maurey asserts that in addition to (\ref{cotype}), we
have a uniform equivalence
$$
\Bignorm{\sum_k\epsilon_k\otimes x_k}_{{\rm
Rad}(X)}\,\asymp\,\biggnorm{\Bigl(\sum_k\vert
x_k\vert^{2}\Bigr)^{\frac{1}{2}}}
$$
for finite families $(x_k)_k$ of $X$ (see e.g. \cite[Thm.
16.18]{DJT}). Thus a bounded linear mapping $u\colon C(K)\to B(X)$
is matricially $R$-bounded if there is a constant $C\geq 0$ such
that for any $n\geq 1$, for any matrix $[f_{ij}]\in M_n(C(K))$ and
for any $x_1,\ldots,x_n\in X$, we have
$$
\biggnorm{\Bigl(\sum_i\Bigl\vert\sum_j
u(f_{ij})x_j\Bigr\vert^2\Bigr)^{\frac{1}{2}}}\,\leq
C\,\bignorm{[f_{ij}]}_{C(K;M_n)}\,\biggnorm{\Bigl(\sum_j\vert
x_j\vert^{2}\Bigr)^{\frac{1}{2}}}.
$$
Mappings satisfying this property were introduced by Simard in
\cite{Sim} under the name of $\ell^2$-cb maps. In this section we
will apply a factorization property of $\ell^2$-cb maps
established in \cite{Sim}, in the case when $X$ is merely an
$L^p$-space.

\bigskip
Throughout this section, we let $(\Omega,\mu)$ be a
$\sigma$-finite measure space. By definition, a density on that
space is a measurable function $g\colon\Omega \to (0,\infty)$ such
that $\|g\|_1 = 1$. For any such function and any $1\leq
p<\infty$, we consider the linear mapping
$$
\phi_{p,g} \colon L^p(\Omega,\mu) \longrightarrow L^p(\Omega,g
d\mu),\qquad \phi_{p,g}(h)=g^{-1/p}h,
$$
which is an isometric isomorphism. Note that $(\Omega,gd\mu)$ is a
probability space. Passing from $(\Omega,\mu)$ to $(\Omega,g
d\mu)$ by means of the maps $\phi_{p,g}$ is usually called a
change of density. A classical theorem of Johnson-Jones \cite{JJ}
asserts that for any bounded operator $T\colon L^p(\mu)\to
L^p(\mu)$, there exists a density $g$ on $\Omega$ such that
$\phi_{p,g}\circ T\circ\phi_{p,g}^{-1}\colon L^p(gd\mu)\to
L^p(gd\mu)$ extends to a bounded operator $L^{2}(gd\mu)\to
L^{2}(gd\mu)$. The next statement is an analog of that result for
$C(K)$-representations.

\begin{prop}\label{Prop Change of density}
Let  $1\leq p < \infty$ and let $u \colon C(K) \to B(L^p(\mu))$ be
a bounded homomorphism. Then there exists a density
$g\colon\Omega\to (0,\infty)$ and a bounded homomorphism $w \colon
C(K) \to B(L^{2}(gd\mu))$ such that
$$
\phi_{p,g} \circ u(f) \circ \phi_{p,g}^{-1} = w(f),\qquad   f\in
C(K),
$$
where equality holds on $L^{2}(gd\mu)\cap L^p(gd\mu)$.
\end{prop}

\begin{proof}
Since $X = L^p(\mu)$ has property $(\alpha)$, the mapping $u$ is
matricially $R$-bounded by Proposition \ref{Prop u matricially
R-bounded}. According to the above discussion, this means that $u$
is $\ell^2$-cb  in the sense of \cite[Def. 2]{Sim}. The result
therefore  follows from \cite[Thms. 3.4 and 3.6]{Sim}.
\end{proof}

We will now focus on Schauder bases on separable $L^p$-spaces. We
refer to \cite[Chap. 1]{LTz1} for general information on this
topic. We simply recall that a sequence $(e_k)_{k\geq 1}$ in a
Banach space $X$ is a basis if for every $x\in X,$ there exists a
unique scalar sequence $(a_k)_{k\geq 1}$ such that $\sum_k a_k
e_k$ converges to $x.$ A basis $(e_k)_{k\geq 1}$ is called
unconditional if this convergence is unconditional for all $x\in
X$. We record the following standard characterization.

\begin{lem}\label{Lem Characterization Unconditional Basis}
A sequence $(e_k)_{k\geq 1}\subset X$ of non-zero vectors is an
unconditional basis of $X$ if and only if $X=\overline{{\rm
Span}}\{e_k\, :\, k\geq 1\}$ and there exists a constant $C\geq 1$
such that for any bounded scalar sequence $(\lambda_k)_{k\geq 1}$
and for any finite scalar sequence $(a_k)_{k\geq 1}$,
\begin{equation}\label{Equ Unconditional 2}
\Bignorm{\sum_{k} \lambda_k a_k e_k} \,\leq C
\sup_k\vert\lambda_k\vert\, \Bignorm{\sum_{k} a_k e_k}.
\end{equation}
\end{lem}

We will need the following elementary lemma.

\begin{lem}\label{Lem Technical Unconditional Basis}
Let $(\Omega,\nu)$ be a $\sigma$-finite measure space, let $1\leq
p < \infty$ and let $Q\colon L^p(\nu)\to L^p(\nu)$ be a finite
rank bounded operator such that $Q_{|L^{2}(\nu)\cap L^p(\nu)}$
extends to a bounded operator $L^{2}(\nu)\to L^{2}(\nu)$. Then
$Q(L^p(\nu)) \subset L^{2}(\nu)$.
\end{lem}

\begin{proof}
Let $E=Q(L^p(\nu)\cap L^{2}(\nu))$. By assumption, $E$ is a finite
dimensional subspace of $L^p(\nu)\cap L^{2}(\nu)$. Since $E$ is
automatically closed under the $L^p$-norm and $Q$ is continuous,
we obtain that $Q(L^p(\nu))=E$.
\end{proof}

\begin{thm}\label{Thm Unconditional basis}  Let $1 \leq p < \infty$
and assume that $(e_k)_{k\geq 1}$ is an unconditional basis of
$L^p(\Omega,\mu).$ Then there exists a density $g$ on $\Omega$
such that $\phi_{p,g}(e_k)\in L^{2}(gd\mu)$ for any $k\geq 1$, and
the sequence $(\phi_{p,g}(e_k))_{k\geq 1}$ is an unconditional
basis of $L^{2}(gd\mu).$
\end{thm}

\begin{proof}
Property (\ref{Equ Unconditional 2}) implies that for any
$\lambda=(\lambda_k)_{k\geq 1}\in\ell^{\infty}$, there exists a
(necessarily unique) bounded operator $T_\lambda\colon L^p(\mu)
\to L^p(\mu)$ such that $T_\lambda(e_k)=\lambda_k e_k$ for any
$k\geq 1$. Moreover $\norm{T_\lambda}\leq C\norm{\lambda}_\infty$.
We can therefore consider the mapping
$$
u\colon\ell^{\infty}\longrightarrow B(L^p(\mu)),\qquad
u(\lambda)=T_\lambda,
$$
and $u$ is a bounded homomorphism. By Proposition \ref{Prop Change
of density}, there is a constant $C_{1}>0$, and a density $g$ on
$\Omega$ such that with $\phi=\phi_{p,g}$, the mapping
$$
\phi T_\lambda\phi^{-1}\colon L^{p}(gd\mu)\longrightarrow
L^{p}(gd\mu)
$$
extends to a bounded operator
$$
S_\lambda\colon L^{2}(gd\mu)\longrightarrow L^{2}(gd\mu)
$$
for any $\lambda\in\ell^{\infty}$, with $\norm{S_\lambda} \leq
C_{1}\norm{\lambda}_{\infty}$.

Assume first that $p\geq 2$, so that $L^{p}(gd\mu)\subset
L^{2}(gd\mu)$. Let $\lambda=(\lambda_k)_{k\geq 1}$ in
$\ell^{\infty}$ and let $(a_k)_{k\geq 1}$ be a finite scalar
sequence. Then $S_\lambda(\phi(e_k)) = \phi
T_\lambda\phi^{-1}(\phi(e_k))= \lambda_{k} \phi(e_{k})$ for any
$k\geq 1$, hence
$$
\Bignorm{\sum_k \lambda_k a_k\phi(e_k)}_{L^{2}(gd\mu)} \,=\,
\Bignorm{S_\lambda \Bigl(\sum_k
a_k\phi(e_k)\Bigr)}_{L^{2}(gd\mu)}\,\leq
C_{1}\|\lambda\|_\infty\,\Bignorm{\sum_k
a_k\phi(e_k)}_{L^{2}(gd\mu)}.
$$
Moreover the linear span of the $\phi(e_k)$'s is dense in
$L^{p}(gd\mu)$, hence in $L^{2}(gd\mu)$. By Lemma \ref{Lem
Characterization Unconditional Basis}, this shows that
$(\phi(e_k))_{k\geq 1}$ is an unconditional basis of
$L^{2}(gd\mu)$.

Assume now that $1\leq p<2$. For any $n\geq 1$, let
$f_n\in\ell^\infty$ be defined by $(f_n)_k=\delta_{n,k}$ for any
$k\geq 1$, and let $Q_n\colon L^p(gd\mu)\to L^p(gd\mu)$ be the
projection defined by
$$
Q_n\Bigl(\sum_k a_k\phi(e_k)\Bigr)\,=\, a_n\phi(e_n).
$$
Then $Q_n=\phi T_{f_n}\phi^{-1}$ hence $Q_n$ extends to an $L^2$
operator. Therefore, $\phi(e_n)$ belongs to $L^{2}(gd\mu)$ by
Lemma \ref{Lem Technical Unconditional Basis}.

Let $p'=p/(p-1)$ be the conjugate number of $p$, let
$(e^{*}_{k})_{k\geq 1}$ be the biorthogonal system of $(e_k)_{k
\geq 1}$, and let $\phi'=\phi^{*-1}$. (It is easy to check that
$\phi'=\phi_{p',g}$, but we will not use this point.) The linear
span of the $e_k^{*}$'s is $w^*$-dense in $L^{p'}(\mu)$.
Equivalently, the linear span of the $\phi'(e_k^{*})$'s is
$w^*$-dense in $L^{p'}(gd\mu)$, hence it is dense in
$L^{2}(gd\mu)$. Moreover for any $\lambda\in\ell^{\infty}$ and for
any $k\geq 1$, we have
$T_{\lambda}^{*}(e_{k}^{*})=\lambda_{k}e_{k}^{*}$. Thus for any
finite scalar sequence $(a_k)_{k\geq 1}$, we have
$$
\sum_{k}\lambda_{k}a_{k}\phi'(e_{k}^{*})\, =\, (\phi
T_\lambda\phi^{-1})^{*}\Bigl(\sum_{k}
a_{k}\phi'(e_{k}^{*})\Bigr)\, =\, S_\lambda^{*} \Bigl(\sum_{k}
a_{k}\phi'(e_{k}^{*})\Bigr).
$$
Hence
$$
\Bignorm{\sum_{k}\lambda_{k}a_{k}\phi'(e_{k}^{*})}_{L^{2}(gd\mu)}\,\leq\,
C_{1}\Bignorm{\sum_{k} a_{k}\phi'(e_{k}^{*})}_{L^{2}(gd\mu)}.
$$
According to Lemma \ref{Lem Characterization Unconditional Basis},
this shows that $(\phi'(e_k^{*}))_{k\geq 1}$ is an unconditional
basis of $L^{2}(gd\mu)$. It is plain that $(\phi(e_k))_{k\geq
1}\subset L^{2}(gd\mu)$ is the biorthogonal system of
$(\phi'(e_k^*))_{k\geq 1}\subset L^{2}(gd\mu)$. This shows that in
turn, $(\phi(e_k))_{k\geq 1}$ is an unconditional basis of
$L^{2}(gd\mu)$.
\end{proof}

We will now establish a variant of Theorem \ref{Thm Unconditional
basis} for non unconditional bases. Recall that if $(e_k)_{n\geq
1}$ is a basis on some Banach space $X$, the projections
$P_N\colon X\to X$ defined by
$$
P_N\Bigl(\sum_k a_k e_k\Bigr)\,=\,\sum_{k=1}^N a_k e_k
$$
are uniformly bounded. We will say that $(e_k)_{k\geq 1} $ is an
$R$-basis if the set $\{P_N \, :\, N\geq 1\}$ is actually
$R$-bounded. It follows from \cite[Thm. 3.9]{CPSW} that any
unconditional basis on $L^p$ is an $R$-basis. See Remark
\ref{Final} (2) for more on this.

\begin{prop}\label{Prop R-basis}
Let $1 \leq p < \infty$ and let $(e_k)_{k\geq 1}$ be an $R$-basis
of $L^p(\Omega,\mu).$ Then there exists a density $g$ on $\Omega$
such that $\phi_{p,g}(e_k)\in L^{2}(gd\mu)$ for any $k\geq 1$, and
the sequence $(\phi_{p,g}(e_k))_{k\geq 1}$ is a basis of
$L^{2}(gd\mu).$
\end{prop}

\begin{proof}
According to \cite[Thm. 2.1]{LMS}, there exists a constant $C\geq
1$ and a density $g$ on $\Omega$ such that with $\phi =
\phi_{p,g}$, we have
$$
\|\phi P_N \phi^{-1}h \|_{2} \leq C \|h\|_{2},\qquad N\geq 1,\
h\in L^{2}(gd\mu)\cap L^{p}(gd\mu).
$$
Then the proof is similar to the one of Theorem \ref{Thm
Unconditional basis}, using \cite[Prop. 1.a.3]{LTz1} instead of
Lemma \ref{Lem Characterization Unconditional Basis}. We skip the
details.
\end{proof}

\begin{rem}\label{Final}\

(1) Theorem \ref{Thm Unconditional basis} and Proposition
\ref{Prop R-basis} can be easily extended to finite dimensional
Schauder decompositions. We refer to \cite[Sect. 1.g]{LTz1} for
general information on this notion. Given a Schauder decomposition
$(X_k)_{k\geq 1}$ of a Banach space $X$, let $P_N$ be the
associated projections, namely for any $N\geq 1$, $P_N\colon X\to
X$ is the bounded projection onto $X_1\oplus\cdots\oplus X_N$
vanishing on $X_k$ for any $k\geq N+1$. We say that $(X_k)_{k\geq
1}$ is an $R$-Schauder decomposition if the set $\{P_N\, :\, N\geq
1\}$ is $R$-bounded. Then we obtain that for any $1<p<\infty$ and
for any finite dimensional $R$-Schauder (resp. unconditional)
decomposition $(X_k)_{k\geq 1}$ of $L^{p}(\mu)$, there exists a
density $g$ on $\Omega$ such that $\phi_{p,g}(X_k)\subset
L^{2}(gd\mu)$ for any $k\geq 1$, and $(\phi_{p,g}(X_k))_{k\geq 1}$
is a Schauder (resp. unconditional) decomposition of
$L^{2}(gd\mu)$.

\smallskip
(2) The concept of $R$-Schauder decompositions goes back (at
least) to \cite{CPSW} and to various works on $L^p$-maximal
regularity and $H^{\infty}$-calculus, see in particular
\cite{KL,KaW}. Let $C_p$ denote the Schatten spaces. For
$1<p\not=2<\infty$, an explicit example of a Schauder
decomposition on $L^2([0,1];C_p)$ which is not $R$-Schauder is
given in \cite[Sect. 5]{CPSW}. More generally, it follows from
\cite{KL} that whenever a reflexive Banach space $X$ has an
unconditional basis and is not isomorphic to $\ell^{2}$, then $X$
admits a finite dimensional Schauder decomposition which is not
$R$-Schauder. This applies in particular to $X=L^p([0,1])$, for
any $1<p\not=2<\infty$. However the question whether $L^p([0,1])$
admits a Schauder basis which is not $R$-Schauder, is apparently
an open question.

We finally mention that according to \cite[Thm. 3.3]{KaW}, any
unconditional decomposition on a Banach space $X$ with property
($\Delta$) is an $R$-Schauder decomposition.
\end{rem}

\bigskip
\

\end{document}